\documentclass[11pt]{amsart}   
\usepackage{amssymb,xy}
\usepackage{amsthm, amsmath}
\usepackage[latin1]{inputenc}
\usepackage[T1]{fontenc}
\xyoption{all}
\theoremstyle{plain}

\theoremstyle{definition}

\theoremstyle{remark}

\usepackage{graphics}

\usepackage{longtable}

\usepackage{color}

\usepackage{epsfig}

\usepackage{amssymb}

\def\qed{\hfill  \framebox(5,5){}}

\def\resultant{{\rm Res}}

\def\cP{{\mathcal P}}
\def\cc{{\mathcal C}}
\def\bb{{\mathcal B}}

\def\cluster{{\mathfrak Cluster}}

\def\qed{\hfill  \framebox(5,5){}}

\def\cPP{\overline{\cP}(t)}
\def\ccc{\overline{\cc}}

\def\R{{\mathbb R}}
\def\C{{\mathbb C}}

\newcommand{\cF}{\mathcal{F}}

\def\sup{{\rm sup}}
\def\inf{{\rm inf}}

\def\lim{{\rm lim}}

\title{Parametrization of $\epsilon$-rational curves: error analysis}
\author{
Sonia L. Rueda
\and
Juana Sendra
}
\date{January 2010}

\address{\small \rm  Dpto. de Matem\'atica Aplicada,
E.T.S. Arquitectura, Univ.
Polit\'ecnica de Madrid, Spain }
\email{sonialuisa.rueda@upm.es}

\address{\small \rm Dpto. de Matem\'atica Aplicada
a la I.T. de Telecomunicaci\'on,
E.U.I.T.Telecomunicaci\'on,
Univ. Polit\'ecnica de Madrid,
Spain }
\email{jsendra@euitt.upm.es}
\usepackage{graphics}

\usepackage{longtable}

\begin{document}

\begin{abstract}
In \cite{PDRuSeSe09} the authors present an algorithm to parametrize approximately $\epsilon$-rational curves, and they show in $2$ examples that the
Hausdorff  distance, w.r.t. to the Euclidean distance, between the input and output curves is small. In this paper, we analyze this distance for a whole family of curves randomly generated and we automatize the strategy used in \cite{PDRuSeSe09}. We find a reasonable upper bound of the Hausdorff distance between each input and output curve of the family.
\end{abstract}

\maketitle

\section*{Introduction}


\noindent
The problem of the approximate parametrization of algebraic plane
curves goes as follows: given a plane algebraic curve $\cc$
(that is the perturbation of a rational plane curve) and a
tolerance $\epsilon>0$, we want to find a new curve
$\overline{\cc}$, being rational, as well as a rational
parametrization of it such that $\cc $ and $\overline{\cc }$ are at certain small distance dependent on
$\epsilon$.

In \cite{PDSeSe04} and \cite{PDSeSe05} it was seen how to approximately
parametrize algebraic plane curves and algebraic surfaces having
an $\epsilon$-singularity of maximum multiplicity. In  \cite{PDRuSeSe09},
using the techniques of $\epsilon$-points developed in \cite{PDSeSe04} and \cite{PDSeSe06}, we extended the results
in \cite{PDSeSe04} to the general case of algebraic affine plane
curves. More precisely, we provide in \cite{PDRuSeSe09} the approximate parametrization
algorithm which given $\cc$ returns an approximate parametrization
$\cPP$ of the curve $\overline{\cc}$.

A natural question arises, the  closeness between the input and output curves of the
algorithm. In our case, this  closeness notion  is given by the
Hausdorff distance. That is, we say that the input and output
curves are close if their Hausdorff distance (as real curves) is small related to the tolerance. We proved in
\cite{PDRuSeSe09} that the Hausdorff distance between $\cc$ and
$\overline{\cc}$ is finite.

In this paper, apply the approximate parametrization algorithm to a
family of curves $\cF$ of degree $4$ which was randomly generated. For each curve
$\cc$ in $\cF$ we compute a bound ${\bb}$ of the Hausdorff distance between $\cc$ and
$\overline{\cc}$. For all of the curves we get ${\bb}\leq 2.2$ and we obtain evidences that the
actual distance is experimentally $\leq 0.1$.

The paper is organized as follows. We recall the approximate parametrization algorithm in
Section \ref{sec-algoritmo}. In Section \ref{sec-familia} we explain how the family ${\mathcal F}$ of curves
of degree $4$ was generated. The last section is devoted to the analysis of the distance between the
the curves of $\cF$ and their approximate parametrizations output by our algorithm.

The following  terminology will be used throughout the paper.
$\|\cdot\|$  and $\|\cdot \|_2$ denote the polynomial  $\infty$--norm and the usual unitary norm in ${\Bbb C}^2$, respectively. $|\cdot|$ denotes the module in $\Bbb C$.
The partial derivatives of a polynomial $g\in {\Bbb C}[x,
y]$  are denoted by
$g^{\overrightarrow{v}}:=\partial^{i+j} g/\partial^{i} x
\partial^{j} y$,
where $\overrightarrow{v}=(i,j)\in {\Bbb N}^{2}$;
we assume that $g^{\tiny{\overrightarrow{0}}}=g.$ Moreover, for
$\overrightarrow{v}=(i,j)\in {\Bbb N}^{2}$,
$|\overrightarrow{v}|=i+j$. Also,   $\overrightarrow{e_1}=(1,0)$, and $\overrightarrow{e_2}=(0,1)$.

\section{Recalling the parametrization algorithm}\label{sec-algoritmo}
In this section we recall the algorithm presented in \cite{PDRuSeSe09} as well as its main properties; see \cite{PDRuSeSe09} for further details. We start with a fixed tolerance $\epsilon$, and with the implicit equation $f(x,y)$ of a real  plane algebraic curve $\cc$ of exact degree $d$, which is the perturbation of a rational curve.  $\cc$ is supposed to satisfy that:
\begin{itemize}
\item[(1)] the degree $d$ of $f$ is proper; this means that there exists a partial derivative of $f$, of order $d$, that  in module is strictly bigger than $\epsilon \|f\|.$
\item[(2)] $f$ is $\epsilon$-irreducible.
\item[(3)] $\cc$ has $d$ different points at infinity, and it does not pass through
$(1 : 0 : 0), (0 : 1 : 0)$.
\end{itemize}
Conditions (1) and (2) guarantee that, under the tolerance, we really have an irreducible curve of degree $d$. Condition (3) ensures that $\cc$ is either compact (as a subset of $\R^2$) or it follows real asymptotes. Therefore, we are excluding curves having a parabolic behavior. The requirement on   $(1 : 0 : 0), (0 : 1 : 0)$ is technical and it can be
achieved by performing a suitable and orthogonal linear change of coordinates.

The theoretical argumentation of the algorithm is as follows.
First, the notion of exact singularity is replaced by the concept of $\epsilon$-singularity, similarly with the notions of exact multiplicity and $\epsilon$-multiplicity. Here, the first complication appears since the number of $\epsilon$-singularities is bigger than (expected) in the exact case; probably due to the perturbation. In order to deal with this difficulty we associate to each $\epsilon$-singularity a radius, and hence we see it as an Euclidean disk. Next, we introduce an equivalence relation on the set of disks and we define the $\epsilon$-singular clusters as the equivalence classes. Then, we define the $\epsilon$-multiplicity of the cluster as the maximum of the $\epsilon$-multiplicities within the class, and we take as canonical representant of the cluster an $\epsilon$-singularity where the $\epsilon$-multiplicity of the cluster is achieved. In this situation, we say that $\cc$ is {\sf $\epsilon$-rational} if the clusters satisfy the well-known genus formula of the exact case. More precisely, if $\{\cluster_{r_i}(Q_i)\}_{i=1,\ldots,s}$ is the cluster
decomposition ($Q_i$ denotes the canonical representant and $r_i$ the $\epsilon$-multiplicity of the cluster),  we say
that $\mathcal C$ is  {\sf $\epsilon$-(affine) rational} if $$
(d-1)(d-2)-\sum_{i=1}^{s} r_i(r_i-1) =0.$$

Now, let us assume that $\cc$ is  $\epsilon$-rational, an let us see how the approximate parametrization algorithm proceeds. The basic idea is, as in the exact case (see \cite{SWP}), to construct a suitable linear system of curves of degree $d-2$.  More precisely, if
$$\{\cluster_{r_i}(Q_i)\}_{i=1,\ldots,s}, \,\,\, Q_i=(q_{i1}:q_{i2}:1)$$ is the cluster decomposition, we compute  $d-3$ simple
$\epsilon$--points on $\mathcal C$, say $\{P_1,\ldots,P_{d-3}\}$ with $P_i=(p_{i1}:p_{i2}:1)$. Again, we associate to each  $\epsilon$-point an Euclidean disk via a radius, and we apply the equivalence relation. If, somehow, any $P_i$ is identified with another $P_j$ or with a singular cluster, we replace $P_i$ by a new  $\epsilon$-point.  In this situation, we consider the effective divisor
$$D=\sum_{i=1}^{s} r_i Q_i +\sum_{i=1}^{d-3} P_i $$
and the (exact) linear system $\mathcal H$ of curves of degree $(d-2)$  given by $D$. That is, $\mathcal H$ is the linear system of curves of degree $d-2$ having  $Q_i$ as $(r_i-1)$-base points,   and $P_i$ as simple base points. If we were working exactly, all intersection points in  ${\mathcal H} \cap \cc$ would be fixed (namely those points in $D$) with the exception of one point that would provide the parametrization. Indeed, in the exact case, the parametrization would  be
$$ \left(\frac{\resultant_y(H(x,y,1), {f})}{\prod_{i=1}^{s}(x-{q}_{i,1})^{r_i(r_i-1)}
\prod_{i=1}^{d-3}(x-{p}_{i,1})}, \frac{
\resultant_x(H(x,y,1), {f})}{\prod_{i=1}^{s}(y-{q}_{i,2})^{r_i(r_i-1)}\prod_{i=1}^{d-3}(y-{p}_{i,2})}\right),$$
being $H(x,y,z)$ the homogeneous polynomial defining $\mathcal H$. In the approximate case, instead of the exact division above, we take the quotient of the Euclidean division of each numerator by the corresponding denominator.

The  output curve derived from this process has the same structure at infinity as the input curve and the same degree, see \cite{PDRuSeSe09}, Theorem 4.5. These properties will play a fundamental role in
the error analysis (see Section \ref{sec-theoretical-error}).  We outline the algorithm derived from the above ideas.

\begin{itemize}
\item[(1)] {\sf Compute} the singular cluster decomposition
$\{\cluster_{r_i}(Q_i)\}_{i=1,\ldots,s}$; say
$Q_i=(q_{i,1}:q_{i,2}:1)$.
\item[(2)] If
$\sum_{i=1}^{s}r_i(r_i-1)\neq (d-1)(d-2),$  {\sf RETURN} ``{\sf
$\mathcal C$ is not (affine) $\epsilon$-rational}". If $s=1$ one may
apply the algorithm in \cite{PDSeSe04}.
\item[(3)] {\sf Compute}
$(d-3)$  $\epsilon$--simple points $\{P_j\}_{1\leq j\leq
d-3}$ of ${\mathcal C}$. Take the points over $\Bbb R$, or as conjugate
complex points. After each point computation check that it is not
in the cluster of the others (including the clusters of $Q_i$); if
this fails take a new one. Say $P_i=(p_{i,1}:p_{i,2}:1)$.
\item[(4)] {\sf Determine} the
linear system  ${\mathcal H}$  of degree $(d-2)$ given
by the divisor $\sum_{i=1}^{s} r_i Q_i+\sum_{i=1}^{d-3} P_i$.
 Let $H(t,x,y,z)=H_1(x,y,z)+tH_2(x,y,z)$
be its defining polynomial.
 \item[(5)] If [$\gcd(F(x,y,0),H_1(x,y,0))\neq 1$] and
[$\gcd(F(x,y,0),H_2(x,y,0))\neq 1$]  replace $H_2$ by
$H_2+\rho_1x^{d-2}+\rho_2y^{d-2}$, where $\rho_1, \rho_2$ are real
and strictly smaller than $\epsilon$.  Say that
$\gcd(F(x,y,0),H_2(x,y,0))=1$; similarly in the other case.
\item[(6)]
$ {S}_1(x, t)=\resultant_y(H(x,y,1), {f})$ and   $ {S}_2(y,
t)=\resultant_x(H(x,y,1), {f}).$
 \item[(8)] $A_1(x)=\prod_{i=1}^{s}(x-{q}_{i,1})^{r_i(r_i-1)}
\prod_{i=1}^{d-3}(x-{p}_{i,1}),$
\\
$A_2(y)=\prod_{i=1}^{s}(y-{q}_{i,2})^{r_i(r_i-1)}\prod_{i=1}^{d-3}(y-{p}_{i,2})$.
\item[(9)]  For $i=1,2$ {\sf compute}  the   quotient $B_i$  of $
{S}_i$ by $A_i$ w.r.t.  either $x$ or $y$. \item[(10)] If the
content of  $B_1$ w.r.t $x$  or the content of $B_2$ w.r.t. $y$ does
depend on $t$, {\sf RETURN} {\sf ``degenerate case"}.
\item[(11)]  {\sf Determine} the root $\overline{p}_1(t)$ of $B_1$,
as a polynomial in $x$, and  the root $\overline{p}_2(t)$ of $B_2$,
as a polynomial in $y$. \item[(12)] {\sf RETURN}  $\overline{\mathcal P}(t)=(\overline{p}_1(t),\overline{p}_2(t))$.
\end{itemize}

\section{Generating a family of $\epsilon$-rational curves}\label{sec-familia}

In this section, we generate the family of curves that will be used in the error analysis.
 We fix three points
$P_1=(2:0:1)$,$P_2=(0:0:1)$ and $P_3=(1:1:1)$ in
$\mathbb{P}^2({\Bbb C})$ and we consider the linear system of curves
of degree 4 defined by the divisor $2P_1+2P_2+2P_3$. Its defining
polynomial is

\noindent
$G(x,y,z,u_1,\ldots,u_6)=u_{2}y^2z^2+u_{3}y^3z+u_{4}y^4+u_{5}xyz^2-
(2u_{2}+3u_{3}+4u_{4}+\frac{1}{2}u_{5}+2u_{6})xy^2z
+u_{6}xy^3+u_{1}x^2z^2+(-\frac{3}{2}u_{5}+2u_{3}
+4u_{4}+2u_{6}-u_{1})x^2yz+  (u_{2}+u_{3}+\frac{1}{2}u_{5}+\frac{1}{4}u_{1}+u_{4})x^2y^2
-u_{1}x^3z+(\frac{1}{2}u_{5}-u_{3}-2u_{4}-u_{6}+\frac{1}{2}u_{1})x^3y+\frac{1}{4}u_{1}x^4.$

Note that for every specialization of $u_i$ such that $G(x,y,z,u_1,\ldots,u_6)$ is irreducible, we get an (exact) rational curve.

Now, for $j=1,\ldots ,6$ and $i=1,\ldots ,10$ let $r_{ij}$ be a random
integer number in the interval $[0,100]$. We obtain 60 different
polynomials $G_{ij}(x,y,z)$, $j=1,\ldots, 6$, $i=1,\ldots ,10$
 setting
$$u_k=\left\{\begin{array}{lcc}(\frac{r_{ij}}{100})^i&\mbox{ if }&k=j\\ 1&\mbox{ if }&k\neq j\end{array}\right. \,\,\, k=1,\ldots ,6$$
in $G(x,y,z,u_1,\ldots,u_6)$.
Given $i\in\{1,\ldots ,6\}$ and $j\in \{1,\ldots ,10\}$ we obtain a random perturbation $g_{ij}(x,y)\in \mathbb{R}[x,y]$ of $G_{ij}(x,y,1)$ as follows
$$g_{ij}(x,y)= G_{ij}(x,y,1)+ \epsilon \frac{r_1}{100} (x+y)+\epsilon^2 \frac{r_2}{100}(x^2+xy+y^2)+$$
$$ \epsilon^3 \frac{r_3}{100}(x^3+x^2y+xy^2+y^3)$$ where
$r_1,r_2,r_3$ are integer numbers taken randomly in the interval
$[0,100]$ and $\epsilon=\frac{1}{100}$. The polynomials $g_{ij}(x,y)$,
$j=1,\ldots, 6$, $i=1,\ldots ,10$ have proper degree 4 and define
60 curves ${\mathcal C}_{ij}$ verifying $(1:0:0),(0:1:0)\notin {\mathcal C}_{ij}^h$  (${\mathcal C}_{ij}^h$ is the projective closure of ${\mathcal C}_{ij}$) and such that they have 4 different points at infinity. Therefore, each of 60 curves satisfies the hypothesis required in parametrization algorithm.

\begin{center}
\begin{figure}[ht]
\centerline{
\psfig{figure=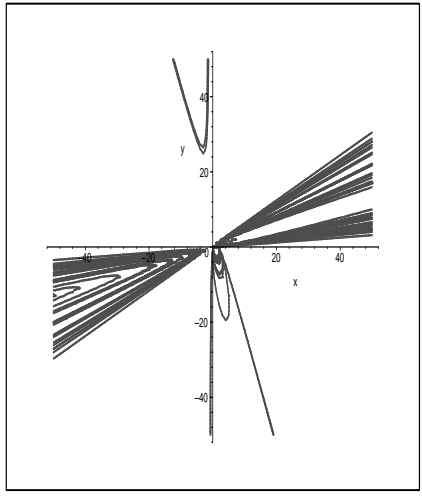,width=6.5cm,height=5.5cm,angle=270}
 }
\caption{ Plot of  the
28 $\epsilon$--rational curves ${\mathcal C}_{ij}$ randomly generated. } \label{fig_fam}
\end{figure}
\end{center}

 Using
the parametrization  algorithm described in Section \ref{sec-algoritmo},  we  conclude that 28 of the 60 curves are
$\epsilon$--rational. We show those curves in Fig. \ref{fig_fam}. An statistical error analysis was given in \cite{PDRuSeSe09b}.
The precise equations of $\cc_{ij}$ as well as the parametrizations provided by the algorithm can be found in

$http://www.aq.upm.es/Departamentos/Matematicas/srueda/fam4.pdf.$

\section{Theoretical strategy for the  error analysis}\label{sec-theoretical-error}

We  describe   the theoretical  strategy for the error analysis that will be employed in the last section. Let $\epsilon>0$ be the tolerance, $\cc$ the input curve and $\ccc$ the output curve provided by the approximate parametrization algorithm. Also, let $f(x,y)$ and $\overline{f}(x,y)$ be the defining polynomials of $\cc$ and $\ccc$, respectively. Moreover, let $\cPP$ be the parametrization of $\ccc$ output by the algorithm. Let $\cc^\R$ and $\ccc^\R$ denote the real part of $\cc$ and $\ccc$, respectively.

Now, we consider the Hausdorff distance $H$, induced by Euclidean distance $d$ in the real plane. That is, for  $A,B\subset \R^2\setminus \{\emptyset\}$
$$ H(A,B)=\max\{\sup_{a\in A}\{d(a,B)\}, \sup_{b\in B}\{d(b,A)\} \}$$
where $d(a,B)=\inf_{b\in B}\{d(a,b)\},  $
and by convection $H(\emptyset,\emptyset)=0$ and, for $\emptyset \neq A\subset X$, $H(A,\emptyset)=\infty$
(see \cite{AB}, for further details).

 The most relevant fact, for our analysis, is that $H(\cc^\R,\ccc^\R)<\infty$ (see \cite{PDRuSeSe09}, Lemma 6.1) because we want to bound $H(\cc^\R,\ccc^\R)$. For this purpose, we will proceed in a similar manner to Section 6 in \cite{PDRuSeSe09}, not only for a couple of examples, but for the whole family of $\epsilon$-rational curves randomly generated in Section \ref{sec-familia}. More precisely, we consider the normal line to $\ccc$ at the generic point $\cPP$, say $ {\mathcal L}_1(t,s)$, as well as the normal line to $\cc$ at the generic point $(a,b)\in \cc^\R$, say  ${\mathcal L}_2(a,b,s)$.
Moreover, we introduce the polynomials
$$ {\mathcal D}_1(t,s)=f({\mathcal L}_1(t,s))\in \overline{\R(t)}[s],\,\,{\mathcal D}_2(a,b,s)=\overline{f}({\mathcal L}_2(a,b,s))\in \C(\overline{\cc})[s],  $$
where $\overline{\R(t)}$ denotes the algebraic closure of $\R(t)$ and $\C(\overline{\cc})$ the field of rational functions over $\overline{\cc}$.
In addition, for every $t_0\in \R$, such that ${\mathcal D}_1(t_0,s)$ is well defined and has real roots,  and  for every $(a_0,b_0)\in \cc^\R$, such that ${\mathcal D}_2(a_0,b_0,s)$ is well defined and has real roots, we take
$$ \begin{array}{l} \rho_{1}^{\R}(t_0)=\min\{|s_0|\, / \, {\mathcal D}_1(t_0,s_0)=0 \,\,\mbox{and}\,\, s_0\in \R\}, \\
 \rho_{2}^{\R}(a_0,b_0)= \min\{|s_0|\, / \, {\mathcal D}_2(a_0,b_0,s_0)=0 \,\,\mbox{and}\,\, s_0\in \R\}. \end{array}  $$
Then,  the supremum of $\rho_{1}^{\R}(t), \rho_{2}^{\R}(a,b)$ provides an upper bound of the Hausdorff distance; at least for those subsets of both curves where the considered minimums are well defined.

Because of computational difficulties, in our analysis, instead of computing $\rho_{1}^{\R}(t), \rho_{2}^{\R}(a,b)$,  we will study
  $$\begin{array}{l} \rho_1(t_0):=\min\{|s_0|\, / \, {\mathcal D}_1(t_0,s_0)=0 \}, \\
  \rho_2(a_0,b_0):=\min\{|s_0|\, / \, {\mathcal D}_2(a_0,b_0,s_0)=0\}. \end{array}
$$
These quantities bound  $d(\overline{\cP}(t_0),\cc)$ and $d((a_0,b_0),\overline{\cc})$ respectively (here, we understand that $d$ is the unitary distance in $\C^2$), instead of $d(\overline{\cP}(t_0),\cc^\R)$ and $d((a_0,b_0),\overline{\cc}^\R)$. So,  for those subsets of both curves, where the corresponding polynomials  are well defined, we bound
$$ \Delta(\cc^\R,\overline{\cc}^{\R}):=\max\{ \sup_{t\in\R}\{d(\cPP,\cc)\}, \sup_{(a,b)\in \overline{\cc}^{\R}}\{d((a,b),\cc)\}\}. $$
For this purpose, in the next section, we will follow the next steps:
\begin{itemize}
\item[(1)] We compute a bound of $\sup_{t\in\R}\{d(\cPP,\cc)\}$, by  applying Corollary 6.2 in \cite{PDRuSeSe09} that ensures that
$$d(\overline{\cP}(t),\cc) \leq  \min\left\{ \left(^{n}_{i}\right)\left| \frac{A_0(t)}{A_i(t)} \right|^{\frac{1}{i}} \,\,\mbox{where}\,\, A_i(t)\neq 0 \,\,\mbox{and}\,\,1\leq i \leq n \right\},$$
where  ${\mathcal D}_1(t,s)=A_n(t) s^n+ \cdots +A_0(t)$.
\item[(2)] In order to bound  $\sup_{(a,b)\in \overline{\cc}^{\R}}\{d((a,b),\cc)\})$ one may apply the same corollary as in (1). However, it implies  to maximize a bivariate rational function under the constrain $f(x,y)=0$. This can be done, for instance, using Lagrange multipliers. Nevertheless, in practice, this is unfeasible. Instead, we estimate the bound by taking a lattice of points $(a,b)\in \ccc$  where we bound $d((a,b),\overline{\cc})$.
\item[(3)]
Note that the quantity  $\Delta(\cc^\R,\overline{\cc}^{\R})$  gives information on how close  every real point on each of the curves is  of a complex point on the other curve. However,   $\Delta(\cc^\R,\overline{\cc}^{\R})\leq H(\cc^\R,\overline{\cc}^{\R})$. To avoid this difficulty, in our analysis, we will look for empirical evidences indicating that the computed bound of $\Delta(\cc^\R,\overline{\cc}^{\R})$ also bounds
$H(\cc^\R,\overline{\cc}^{\R})$; for that we test empirically that, in our computations, $\rho_1(t_0)=\rho_{1}^{\R}(t_0), \rho_{2}(a_0,b_0)=\rho_{2}^{\R}(a_0,b_0)$.

\item[(4)] It may happen for some $t_0\in\R$ that $\rho_{1}^{\R}(t_0)$ is not well defined or simply that the bound it provides of $d(\overline{\cP}(t_0),\cc^\R)$ is not satisfactory.
Observe that to bound $d(\overline{\cP}(t_0),\cc^\R)$ we can use the intersection of any line through $\overline{\cP}(t_0)$ with $\cc$. So in some cases
we will also proceed in the following way.

Let us consider the line at the generic point $\overline{\cP}(t)$ in the direction given by $h\in\R$
\[ {\mathcal L}_{h}(t,s)=\left(\overline{p_1}(t)+s \frac{2h}{h^2+1},\overline{p_2}(t)+s \frac{h^2-1}{h^2+1}\right). \]
We introduce the polynomial ${\mathcal D}_{h}(t,s)=f({\mathcal L}_{h}(t,s))\in \overline{\R(t,h)}[s]$.
For a fixed $h_0\in \R$ and for every $t_0\in \R$, such that ${\mathcal D}_{h_0}(t_0,s)$ is well defined and has real roots,
$d(\overline{\cP}(t_0),\cc^\R)\leq  \rho_{h_0}^{\R}(t_0)$, where
\[\rho_{h_0}^{\R}(t_0)=\min\{|s_0|\, / \, {\mathcal D}_{h_0}(t_0,s_0)=0 \,\,\mbox{and}\,\, s_0\in \R\}.\]
Thus,  the supremum of $\rho_{h_0}^{\R}(t)$ and $\rho_{2}^{\R}(a,b)$ provides an upper bound of the Hausdorff distance.
Then for $t_0\in \R$ such that ${\mathcal D}_{h_0}(t_0,s)$ is well defined, we may study
\[  \rho_{h_0}(t_0):=\min\{|s_0|\, / \, {\mathcal D}_{h_0}(t_0,s_0)=0 \} \]
which is an upper bound of $d(\overline{\cP}(t_0),\cc)$ and plays the role of $\rho_{1}(t_0)$ in the previous steps.
\end{itemize}

\section{Execution of the error analysis}

Let $\cF$ be the family of $28$ $\epsilon$-rational curves of degree $4$ defined in Section 3.
For each curve $\cc$ in $\cF$ we explain next how the bound of $\Delta(\cc^\R,\overline{\cc}^{\R})$ was computed  and show evidences that this bound is also an upper bound of the Hausdorff distance $H(\cc^\R,\overline{\cc}^{\R})$.

\subsection{Bound of $\sup_{t\in\R}\{d(\cPP,\cc)\}$}\label{subsec_bound}

Let us denote by $\Lambda$ the domain of ${\mathcal D}_1(t_0,s)$. If the curve is
compact then the polynomial ${\mathcal D}_1(t_0,s)$ is well defined for every $t_0\in \R$, otherwise it is not defined for
two real poles $\beta_1, \beta_2$ of $\overline{\cP}(t)$. In the family $\cF$ only two curves are compact.

In order to bound $\sup_{t\in\R}\{d(\cPP,\cc)\}$, we obtain an upper  bound of $\rho_1(t)$
when $t\in \Lambda$. For this purpose we maximize the functions $R_1(t)=4\left| \frac{A_0(t)}{A_1(t)} \right|$ and $R_2(t)=\left(^{4}_{2}\right)\left| \frac{A_0(t)}{A_2(t)} \right|^{\frac{1}{2}}$ (see step (1) in Section \ref{sec-theoretical-error}) in $\Lambda$ as
follows.

Let $\alpha_1$ and $\alpha_2$ be the  real roots of the denominator of
$R_1(t)$. For all the curves in $\cF$, $\alpha_1$ and $\alpha_2$ are not real roots
of the denominator of $R_2(t)$.
Let  $I_i$, $i=1,2$ be an interval isolating $\alpha_i$ from $\alpha_j$, $j\neq i$
and from the real poles of $R_2(t)$.  Observe that
$R_1(t)$ and $R_2(t)$ are continuous in ${\Bbb R}\setminus
(I_1\cup I_2)$ and in the adherence $\overline{I_1\cup I_2}$ of $I_1\cup I_2$,
respectively. We compute
$${\mathcal B}_1= \max\{R_1(t)\,|\, t\in {\Bbb R}\setminus (I_1\cup I_2) \}\,\,\mbox{ and }\,\,{\mathcal B}_2=\max\{R_2(t)\,|\, t\in \overline{I_1\cup I_2} \}.$$
Then ${\mathcal B}=\max \{{\mathcal B}_1,{\mathcal B}_2\}$ is an upper bound of $\rho_1(t)$ in $\Lambda$.

The last column of the next table contains the computed bound ${\mathcal B}$ for each one of the curves $\cc_i$, $i=1,\ldots ,28$ of the family $\cF$.
\begin{center}
 Table 1:
\end{center}
\begin{center}
\begin{tabular}{|c|c|c|c|}
  \hline

$i$ & ${\mathcal B}_1$ & ${\mathcal B}_2$ & ${\mathcal B}$ \\
 \hline
1 & 0.3012751472 & 1.784885546 & 1.784885546 \\
2 & 0.1680336313 & 0.8228821157 & 0.8228821157 \\
3 & 0.2209183305 & 1.143210796 & 1.143210796 \\
4 & 0.2457462218 & 1.388890611 & 1.388890611 \\
5 & 0.4775061243 & 1.471164469 & 1.471164469 \\
6 & 0.1854050321 & 0.9172323537 & 0.9172323537 \\
7 & 0.3392516285 & 1.238494405 & 1.238494405 \\
8 & 0.1687631697 & 0.9278483955 & 0.9278483955 \\
9 & 0.4481254299 & 1.345341665 & 1.345341665 \\
10 & 0.1706747632 & 1.252669418 & 1.252669418 \\
11 & 0.4336254993 & 1.328637472 & 1.328637472 \\
12 & 0.04502452088 & 0.7028506083 & 0.7028506083 \\
13 & 0.2511290220 & 1.849173820 & 1.849173820 \\
14 & 0.7973544750 & 0.5426224779 & 0.7973544750 \\
15 & 1.947190823 & 1.201605769 & 1.947190823 \\
16 & 0.1658993167 & 2.124343900 & 2.124343900 \\
17 & 0.06346428265 & 1.634020447 & 1.634020447 \\
18 & 1.401107905 & 1.830395156 & 1.830395156 \\
19 & 1.690902532 & 0.9731554792 & 1.690902532 \\
  20 & 1.590285558 & 1.659392056 & 1.659392056 \\
21 & 0.1543267485 & 0.9761129297 & 0.9761129297 \\
22 & 0.8247063503 & 0.8459442935 & 0.8459442935 \\
23 & 0.2452763324 & 0.6140170288 & 0.6140170288 \\
24 & 0.08434729326 & 0.7159251709 & 0.7159251709 \\
25 & 0.6464253153 & 2.150679036 & 2.150679036 \\
\hline
\end{tabular}
\end{center}
\newpage
\begin{center}
 Table 1: continued
\end{center}
\begin{center}
\begin{tabular}{|c|c|c|c|}
  \hline
26 & 3.604620794 & 1.418101314 & 3.604620794 \\
27 & 0.08779082555 & 0.7809344831 & 0.7809344831 \\
28 & 0.3673519642 & 1.773877016 & 1.773877016\\
  \hline
\end{tabular}
\end{center}

 We will improve next the bound given for curves $\cc_{16}$ and
$\cc_{26}$. For a fixed $h_0\in \R$ we can write ${\mathcal D}_{h_0}(t,s)=B_n(t) s^n+ \cdots +B_0(t)$ to which Corollary 6.3 in
\cite{PDRuSeSe09} applies. Hence we can obtain an upper bound ${\mathcal B}^{h_0}$ of $\sup_{t\in\R}\{d(\cPP,\cc)\}$ maximizing the new
functions
\[4\left| \frac{B_0(t)}{B_1(t)} \right|\,\,\mbox{ and } \left(^{4}_{2}\right) \left| \frac{B_0(t)}{B_2(t)}\right|^{\frac{1}{2}}\]
in the domain of ${\mathcal D}_{h_0}(t,s)$. As described earlier for $R_1(t)$ and $R_2(t)$ we obtain respectively ${\mathcal B}_1^{h_0}$ and ${\mathcal B}_2^{h_0}$.

For curves $\cc_{16}$ and $\cc_{26}$ we computed the upper bound ${\mathcal B}^{h_0}$ of $\sup_{t\in\R}\{d(\cPP,\cc)\}$ using different values of $h_0$ and we found bounds improving the ones given earlier for the values of $h_0$ shown in the next table.
\begin{center}
 Table 2:
\end{center}
\begin{center}
\begin{tabular}{|c|c|c|c|c|}
  \hline
$i$ & $h_0$ &${\mathcal B}_1^{h_0}$ &${\mathcal B}_2^{h_0}$ & ${\mathcal B}^{h_0}$ \\
  \hline
16 & -1 & 0.5870746534 & 1.287063889 &  1.287063889 \\
26 & $\frac{1}{20}$ & 0.2525792337 &  0.001139009266 & 1.178706930 \\
\hline
\end{tabular}
\end{center}

\subsection{Empirical bound of $\sup_{(a,b)\in \overline{\cc}^{\R}}\{d((a,b),\cc)\})$}

In this section we estimate the bound of $\sup_{(a,b)\in \overline{\cc}^{\R}}\{d((a,b),\cc)\})$.
We estimate the bound by taking a lattice of points $(a,b)\in \ccc$
where we bound $d((a,b),\overline{\cc})$ estimating $\rho_{2}^{\R}(a,b)$.
We show evidences for $\rho_{2}^{\R}(a,b)$ being small and for
$\rho_2(a,b)=\rho_{2}^{\R}(a,b)$.

If the curve is not compact, first we analyze the behavior of the input and output curves
through the real asymptotes. Let us suppose that $\cc$ is a non compact curve in $\cF$ and let ${\mathcal L}_1$ and ${\mathcal L}_2$
be its real asymptotes. By \cite{PDRuSeSe09}, Corollary 4.6 the real asymptotes of $\cc$ and $\overline{\cc}$ are parallel lines
so the Hausdorff distance between them can be easily computed. Let $\overline{\mathcal L}_1$ and $\overline{\mathcal L}_2$  be
the real asymptotes of  $\overline{\cc}$ parallel to ${\mathcal L}_1$ and ${\mathcal L}_2$ respectively. We the value of
\[\eta=\max \{H({\mathcal L}_1,\overline{\mathcal L}_1), H({\mathcal L}_2,\overline{\mathcal L}_2))\}\]
for all the non compact curves of ${\mathcal F}$ in the next table. Then we proceed as follows:
\begin{itemize}
\item[(1)] For each  negative integer  $i$ we compute the set $\Omega_i$ of intersections of
$\cc^\R$ with the line $x=i$. We obtain
$m_{i}^{\R}:=\max\{\rho_{2}^{\R}(a,b)\,|\, (a,b)\in \Omega_i\}$ and
$m_{i}:=\max\{\rho_{2}(a,b)\,|\, (a,b)\in \Omega_i\}$, and we check
that $m_{i}^{\R}=m_{i}$.

\item[(2)] We repeat  the previous step until
\[\min\{|\rho_{2}^{\R}(i,b)- H({\mathcal L}_1,\overline{\mathcal L}_1)|, |\rho_{2}^{\R}(i,b)- H({\mathcal L}_2,\overline{\mathcal L}_2)|\, / \, (i,b)\in\Omega_i\}<\varepsilon.
\]

\item[(3)] Let $\tau_1$ be the smallest  value  of $i$ until termination of this process.
\end{itemize}

We perform this experiment also for each  positive  integer  $i$ to
obtain in this case the highest  value $\tau_2$ such that the inequality in step (2) holds. At the same time we check that $m_{i}^{\R}=m_{i}$ with
positive integers $i=1,\ldots ,\tau_2$. The same process is repeated
for $y=j$,  to obtain the negative and positive integers  $ \tau_3, \tau_4$, respectively such that
\[
\min\{|\rho_{2}^{\R}(a,j)- H({\mathcal L}_1,\overline{\mathcal L}_1)|, |\rho_{2}^{\R}(a,j)- H({\mathcal L}_2,\overline{\mathcal L}_2)|\, / \, (a,j)\in\Omega^j\}<\varepsilon
\]
where $\Omega^j$ is the set of intersections of $\cc^\R$ with the line $y=j$. Let ${(m^{j})}^{\R}:=\max\{\rho_{2}^{\R}(a,b)\,|\, (a,b)\in \Omega^j\}$ and
$m^{j}:=\max\{\rho_{2}(a,b)\,|\, (a,b)\in \Omega^j\}$. We also check that ${(m^{j})}^{\R}=m^{j}$ with integers $j=\tau_3,\ldots,-1,1,\ldots,\tau_4$. Let
$m:=\max\{m_i, m^j \,|\,i = \tau_1,\ldots,-1,1,\ldots,\tau_2, j=\tau_3,\ldots,-1,1,\ldots,\tau_4\}$.

Let $[\tau_1,\tau_2]\times[\tau_3,\tau_4]$. We
empirically consider that out of the compact $[\tau_1,\tau_2]\times[\tau_3,\tau_4]$, the curves behave
as the asymptotes, and  the  empirical bound of $\sup_{(a,b)\in \overline{\cc}^{\R}}\{d((a,b),\cc)\})$
in  $[\tau_1,\tau_2]\times[\tau_3,\tau_4]$ is $m$.

The next table shows the compact set $[\tau_1,\tau_2]\times[\tau_3,\tau_4]$ obtained for $\varepsilon=10^{-6}$ in all the curves
except for numbers $25$ and $28$ for which we took $\varepsilon=10^{-5}$. The reason being that we run out of memory before
reaching the box outside of which the curves behaved like the asymptotes
with $\varepsilon=10^{-5}$.

\begin{center}
 Table 3:
\end{center}
\begin{center}
\begin{tabular}{|c|c|c|c|c|}
  \hline
$i$  & $[\tau_1,\tau_2]$ & $[\tau_3,\tau_4]$ & $m$ & $\eta$\\
\hline
1  & [-3434, \,3428] & [-2069, \,2066] & 0.04474051996 & 0.002685992105 \\
2  & [-6730, \,6732] & [-3833, \,3835] & 0.01909150476 & 0.007250422655\\
3  & [-2120, \,2120] & [-1030, \,1031] & 0.02523781400 & 0.0008098244306\\
4  & [-2485, \,2492] & [-1370, \,1374] & 0.03061351675 & 0.0006391265474\\
5  & [-4157, \,4160] & [-2000, \,2003] & 0.03559861599 & 0.0007910063013\\
6  & [-11115, \,11089] & [-4663, \,4642] & 0.02312094507 & 0.005224665954\\
7  & [-22154, \,22121] & [-9370, \,9334] & 0.03942395261 & 0.01044227346\\
8  & [-11207, \,11214] & [-4756, \,4761] & 0.02118572337 & 0.005399010161\\
9  & [-28777, \,28684] & [-11887, \,11825] & 0.05222449790 & 0.01308589690\\
10  & [-1922, \,1925] & [-1063, \,1065] & 0.02633886946 & 0.0003845822760\\
11  & [-4555, \,4552] & [-1912, \,1912] & 0.03226327617 & 0.001373344655\\
12  & [-5017, \,4997] & [-2114, \,2107] & 0.01246202377 & 0.002415851115\\
14  & [-797, \,780] & [-242, \,235] & 0.01907109331 & 0.001277331144\\
15  & [-4809, \,4774] & [-1410, \,1396] & 0.07529126612 & 0.005987329671\\
16  & [-3841, \,3866] & [-15344, \,15340] & 0.05528214429 & 0.01084165622\\
17  & [-139, \,16] & [-955, \,832] & 0.03702449872 & 0.008212462120\\
19  & [-13558, \,13560] & [-6428, \,6446] & 0.03068560995 & 0.01035036823\\
20  & [-4290, \,4291] & [-2287, \,2288] & 0.04051580912 & 0.001538264208\\
21  & [-14465, \,14520] & [-5275, \,5300] & 0.01835471004 & 0.005065273865\\
22  & [-2286, \,2282] & [-982, \,982] & 0.01911087028 & 0.0002824836230\\
23  & [-946, \,940] & [-400, \,399] & 0.01154559037 & 0.0004413541720\\
24  & [-157, \,3724] & [-1215, \,1213] & 0.01530235590 & 0.0001629862393\\
25  & [-18910,\, 18839]& [-5882, \,5863] & 0.08880952924 & 0.03464418857\\
26  & [-3400, \, 3398] & [1769, \, 1769]& 0.03324362713 & 0.001139009266\\
27  & [-308, \,5274] & [1659, \, 1642]& 0.01648328102 & 0.0002351747177\\
28  & [-6279, \,6210] & [-2001, \, 1982] & 0.03876376237 & 0.01250853150\\
\hline
\end{tabular}
\end{center}
\newpage
If the curve $\cc$ is compact we consider a compact set
$[\tau_1,\tau_2]\times[\tau_3,\tau_4]$ containing  $\cc^{\R}$. Then
we compute $m$ as previously described checking also that
$m_{i}^{\R}=m_{i}$,${(m^{j})}^{\R}=m^{j}$ with $i =
\tau_1,\ldots,-1,1,\ldots,\tau_2$,
$j=\tau_3,\ldots,-1,1,\ldots,\tau_4$.

\vspace*{1 cm}
\begin{center}
 Table 4:
\end{center}
\begin{center}
\begin{tabular}{|c|c|c|c|}
  \hline
$i$ & $[\tau_1,\tau_2]$ & $[\tau_3,\tau_4]$ & $m$ \\
\hline
13  & $[-9/512,8041/1024]$ & $[-20057/1024,2117/1024] $ & 0.04595703645 \\
18  & $[-5/512,592745/512]$ & $[-8723/1024,304847/1024]$ & 0.09228397972\\
\hline
\end{tabular}
\end{center}
\vspace*{1 cm}

\subsection{Empirical evidences}

Now, we perform some empirical tests to show evidences that
$\rho_{1}^{\R}(t)$ is smaller than the upper bound ${\mathcal B}$ of
$\sup_{t\in\R}\{d(\cPP,\cc)\}$ given in Section \ref{subsec_bound}.
First, let $D_{1}(s)=\lim_{t\mapsto \pm \infty} {\mathcal D}_{1}(t,s)$.
Then, for every curve of the family ${\mathcal F}$ let
$$\chi=\min\{|s_0|\, / \, D_1(s_0)=0 \,\,\mbox{and}\,\, s_0\in
\R\}.$$ We checked that $\chi=\min\{|s_0|\, / \, D_1(s_0)=0\}$ in
all cases. Since the roots of a polynomial depend continuously on
its coefficients, for every $\delta>0$ there exists  $K>0$ such that
for all $|t_0|>K$ there is a root $s_0$ of ${\mathcal D}_1(t_0,s)$ with
$\|\chi-s_0\|_2 <\delta$. It may happen that these roots are all
complex. However, in our example, we see that
$\rho_1((-10)^k)=\rho_{1}^{\R}((-10)^k)$ for $k=1,\ldots,20$. Let
$\chi_1=\min\{\rho_{1}^{\R}((-10)^k)\,|\, k=1,\ldots,20\}$ and
$\chi_2=\max\{\rho_{1}^{\R}((-10)^k)\,|\, k=1,\ldots,20\}$. We show
these computations in the next table.

\begin{center}
 Table 5:
\end{center}
\begin{center}
\begin{tabular}{|c|c|c|c|c|}
  \hline
$i$ & $\chi$ & $\chi_1$ & $\chi_2$ \\
\hline
1 & 0.001918863706 & 0.001918568088 & 0.001922644324 \\
2 & 0.004169957700 & 0.004161970065 & 0.004170583666 \\
3 & 0.0006994105148 & 0.0006993543405 & 0.0007001275847 \\
4 & 0.0006662568567 & 0.0006659169716 & 0.0006706114665 \\
5 & 0.001185494963 & 0.001184862374 & 0.001191850748 \\
6 & 0.0002278441391 & 0.0002266218453 & 0.0002434582060 \\
7 & 0.004927221514 & 0.004922425041 & 0.004975309877 \\
8 & 0.001581018062 & 0.001580991493 & 0.001581020714 \\
9 & 0.003893937744 & 0.003759724534 & 0.008421643120 \\
10 & 0.0004038643651 & 0.0004036487366 & 0.0004066244124 \\
11 & 0.001324473283 & 0.001324364325 & 0.001325864354 \\
12 & 0.0006513574356 & 0.0006500882228 & 0.0006514569249 \\
13 & 0.004018853976 & 0.004018487265 & 0.004023539023 \\
14 & 0.001399552125 & 0.001399552100 & 0.001399552395 \\
15 & 0.006540790926 & 0.006540789722 & 0.006540806320 \\
16 & 0.01873734902 & 0.001829670103 & 0.03745550824 \\
17 & 0.0003018645266 & 0.0002989553843 & 0.0003020916294 \\
18 & 0.003558391244 & 0.003557440039 & 0.003570581899 \\

\hline
\end{tabular}
\end{center}
\newpage
\begin{center}
 Table 5: continued
\end{center}
\begin{center}
\begin{tabular}{|c|c|c|c|c|}
  \hline
$i$ & $\chi$ & $\chi_1$ & $\chi_2$ \\
\hline
19 & 0.0005206627618 & 0.0005145010748 & 0.0005996180488 \\
20 & 0.001782196531 & 0.001779734571 & 0.001782429636 \\
21 & 0.001157562963 & 0.001155811228 & 0.001179946145 \\
22 & 0.001096140685 & 0.001095297129 & 0.001096220492 \\
23 & 0.0004670320359 & 0.0004670115956 & 0.0004670339681 \\
24 & 0.001057066460 & 0.001056738383 & 0.001061264741 \\
25 & 0.06749573507 & 0.005148188392 & 0.06814019159 \\
26 & 0.001564435599 & 0.001563775730 & 0.001572903801 \\
27 & 0.001268468809 & 0.001267989563 & 0.001274605947 \\
28 & 0.0005550274941 & 0.0004502065315 & 0.0005632086626\\
\hline
\end{tabular}
\end{center}

 Observe that for most of the curves of the family
$\cF$ the bound ${\mathcal B}$ computed in Section 5.1  was reached at
${\mathcal B}_2$. Let us suppose that the isolating intervals of the
real poles of $R_1(t)$ are $I_1=[t_{11},t_{12}]$ and
$I_2=[t_{21},t_{22}]$. Then if $\Gamma_2=\{t\in\R\,/\, R'_2(t)=0\}$
then
\[{\mathcal B}_2=max\{R_2(t)\,/\, t\in \Gamma_2\cup\{t_{11},t_{12},t_{21},t_{22}\}\}.\]
Furthermore, only for curve $\cc_{25}$ the set $\Gamma_2$ is
nonempty. We show in the next table that
$\nu=\max\{\rho_{1}^{\R}(t_{11}),\rho_{1}^{\R}(t_{12}),\rho_{1}^{\R}(t_{21}),\rho_{1}^{\R}(t_{22})\}$
is very small compared to the value of
$\mu=\max\{R_2(t_{11}),R_2(t_{12}),R_2(t_{21}),R_2(t_{22})\}$ and
also we compare them with ${\mathcal B}_2$.

\begin{center}
 Table  6:
\end{center}
\begin{center}
\begin{tabular}{|c|c|c|c|}
  \hline
$i$ & ${\mathcal B}_2$ & $\mu$ & $\nu$ \\
\hline
1 & 1.784885546 & 1.784885546 & 0.01895037003 \\
2 & 0.8228821157 & 0.8228821157 & 0.003187256654 \\
3 & 1.143210796 & 1.143210796 & 0.007932988669 \\
4 & 1.388890611 & 1.388890611 & 0.01141905470 \\
5 & 1.471164469 & 1.471164469 & 0.01307915412 \\
6 & 0.9172323537 & 0.9172323537 & 0.004427051081 \\
7 & 1.238494405 & 1.238494405 & 0.008019188716 \\
8 & 0.9278483955 & 0.9278483955 & 0.004590148599 \\
9 & 1.345341665 & 1.345341665 & 0.009661551521 \\
10 & 1.252669418 & 1.252669418 & 0.009307873805 \\
11 & 1.328637472 & 1.328637472 & 0.01086359027 \\
12 & 0.7028506083 & 0.7028506083 & 0.002603616572 \\
13 & 1.849173820 & 1.849173820 & 0.01407639891 \\
14 & 0.5426224779 & 0.5426224779 & 0.002434019507 \\
15 & 1.201605769 & 1.201605769 & 0.01147127202 \\
16 & 2.124343900 & 2.124343900 & 0.01663838395 \\
17 & 1.634020447 & 1.634020447 & 0.01149183389 \\
18 & 1.830395156 & 1.830395156 & 0.01586113294 \\
19 & 0.9731554792 & 0.9731554792 & 0.006071538946 \\
20 & 1.659392056 & 1.659392056 & 0.01644172109 \\
21 & 0.9761129297 & 0.9761129297 & 0.005852977427 \\
22 & 0.8459442935 & 0.8459442935 & 0.004712176306 \\
23 & 0.6140170288 & 0.6140170288 & 0.002448022081 \\
 \hline
\end{tabular}
\end{center}

\newpage
\begin{center}
 Table 6: continued
\end{center}
\begin{center}
\begin{tabular}{|c|c|c|c|}
  \hline
$i$ & ${\mathcal B}_2$ & $\mu$ & $\nu$ \\
\hline
24 & 0.7159251709 & 0.7159251709 & 0.006089518478 \\
25 & 2.150679036 & 2.150679036 & 0.05229161121 \\
26 & 1.418101314 & 1.418101314 & 0.01218396938 \\
27 & 0.7809344831 & 0.7809344831 & 0.007553005576 \\
28 & 1.773877016 & 1.773877016 & 0.01971360528\\
 \hline
\end{tabular}
\end{center}

For curves $\cc_i$, $i=14,15,19,26$ the bound ${\mathcal B}$ is achieved
in ${\mathcal B}_1$. Let $\Gamma_1=\{t\in\R\,/\, R'_1(t)=0\}$, then
\[{\mathcal B}_1=\max\{\{R_1(t)\,/\, t\in \Gamma_1\cup\{t_{11},t_{12},t_{21},t_{22}\}\}\cup \{\lim_{t\mapsto \pm \infty} R_1(t)\}\}.\] In fact, in those curves
${\mathcal B}_1$ is equal to $\max\{R_1(t)\,/\, t\in \Gamma_1\}$.

For each real pole of $R_1(t)$ as well as for each real critical
value of $R_1(t)$, we consider a sequence of isolating intervals
$J_k$ of length  $1/10^{k+5}$, we take the middle point $t_k$, and
we analyze $\rho_1(t_k), \rho_{1}^{\R}(t_k)$. After a certain $k_1$
the sequences become stable, let $E_{k_1}$ be the set containing
$\rho_{1}^{\R}(t_{k_1})$ for the $k_1$th element of each one of the
sequences constructed for each real pole of $R_1(t)$. Let
$\gamma_1=\max E_{k_1}$. Similarly, after a certain $k_2$ the
sequences for the real critical values of $R_1(t)$ become stable.
Let $E_{k_2}$ be the set containing $\rho_{1}^{\R}(t_{k_2})$ for the
$k_2$th element of each one of the sequences constructed for each
real critical values of $R_2(t)$. We call $\gamma_2=\max E_{k_2}$.
If the curve is not compact, we also perform this experiment for
each of the two real poles $\beta_1,\beta_2$ of $\overline{\cP}(t)$.
The sequences to $\beta_1,\beta_2$ become stable after a certain
$k_3$ and we get a set $E_{k_3}$ and $\gamma_3=\max E_{k_3}$.

\begin{center}
 Table 7:
\end{center}
\begin{center}
\begin{tabular}{|c|c|c|c|}
  \hline
i & $\gamma_1$ & $\gamma_2$ & $\gamma_3$ \\
\hline
1 & 0.05760637790 & 0.07103885930 & 0.002685991017 \\
2 & 0.03680694646 & 0.04030133537 & 0.007250422354 \\
3 & 0.03985346560 & 0.05219935882 & 0.0008098246652 \\
4 & 0.04358984318 & 0.06256427563 & 0.0006391277723 \\
5 & 0.05222536651 & $\max\{0.009579373594,s_{c_{1}}\}$ & 0.0007910068620 \\
6 & 0.03997031881 & 0.03468777768 & 0.005224667630 \\
7 & 0.05420521510 & 0.06230995426 & 0.01044227515 \\
  8 & 0.03842889832 & 0.03215152707 & 0.005399009187 \\
9 & 0.06425437603 & 0.08158624352 & 0.01308589929 \\
10 & 0.03857067352 & 0.03800183984 & 0.0003845830022 \\
11 & 0.05050891041 & $\max\{0.007217929868, s_{c_5}\}$ & 0.001373344987 \\
12 & 0.02590356741 & 0.01124398891 & 0.002415851766 \\
13 & 0.07620545140 & 0.05039734679 &  \\
14 & 0.02750848667 & $\max\{0.01846235849,s_{c_5}\}$  & 0.001277332295 \\
15 & 0.06016762755 & 0.2971598318 (*) & 0.005987330302 \\
16 & 0.09714515451 & 0.03765565867 & 0.01084165864 \\
 \hline
\end{tabular}
\end{center}

\newpage
\begin{center}
 Table 7: continued
\end{center}
\begin{center}
\begin{tabular}{|c|c|c|c|}
 \hline
17 & 0.06901149454 & 0.01524196447 & 0.008212460400 \\
18 & 0.07011950234 & $\max\{0.01367185130,s_{c_5}\}$ &  \\
19 & 0.05269958334 & 0.2771889722 (*) & 0.01035036672 \\
20 & 0.05744093868 & $\max\{0.01461432871, s_{c_2}\}$ & 0.001538264277 \\
21 & 0.05357323512 & 0.03486848807 & 0.005065270849 \\
22 & 0.04568079074 & $\max\{0.009770701627, s_{c_2}\}$ & 0.0002824834209 \\
23 & 0.02825492360 & 0.06513903382 & 0.0004413527932 \\
24 & 0.05459990275 & 0.01907914217 & 0.0001629875241 \\
25 & 0.09628252969 & 0.1341804080 & 0.03464419053 \\
26 & 0.05149574516 & 1.431046152 (*) & 0.001139009002 \\
27 & 0.06133149851 & 0.01990495615 & 0.0002351734340 \\
28 & 0.09919600166 & 0.09135201774 & 0.01250853076\\
 \hline
\end{tabular}
\end{center}

\subsection{Behaviour of $\rho_{1}^{\R}(t)$}

In the grate majority of our computations $\rho_{1}^{\R}(t)$ is defined and $\rho_{1}^{\R}(t)=\rho_{1}(t)$.
In some cases for a given $t_0\in\R$ the set
$\{|s_0|\, / \, {\mathcal D}_1(t_0,s_0)=0 \,\,\mbox{and}\,\, s_0\in \R\}$
is empty or $\rho_{1}^{\R}(t)$ and $\rho_{1}(t)$ happen to be
different. Then it should  be taken into consideration that there exists $h_0\in\R$
such that $\rho_{h_0}^{\R}(t_0)<\infty$.

Let $s_c$ denote $\{\rho_{1}^{\R}(t_k)\}$ for the sequence $\{t_k\}$
to the critical point $c$ of $R_1(t)$. For curves number $5$, $11$, $18$, $20$ and $22$  the sequence $s_c$ could not
be computed  for some critical point $c$ of $R_1 (t)$. For example, curve $\cc_5$ has $6$ critical points and
the sequence $s_{c_{1}}$ for critical point $c_1$ was not defined. We write $\max\{0.009579373594,s_{c_{1}}\}$ where
$0.009579373594$ is the maximum of the values at which the rest of the sequences stabilized.

The results for $\gamma_2$ marked with (*) indicate that $\rho_{1}^{\R}(t_k)\neq\rho_{1}(t_)$ for the sequence $\{t_k\}$ to one of the critical points of $R_1(t)$. For curve number $26$ there are $6$ critical values of $R_1(t)$ and  $1.357539211\leq |\rho_{1}^{\R}(t_k)-\rho_{1}(t_k)| \leq 1.357539223$ for the sequence $\{t_k\}$ to the critical value $c_5$.

For each one of the curves highlighted and for the critical point of $R_1 (t)$ where $\rho_{1}^{\R}(t)$ did not behave properly (there was only one of those points of each curve) we proceed as follows. We consider a sequence of isolating intervals $J_k$ of length  $1/10^{k+5}$,
we take the middle point $t_k$, and we analyze $\rho_{h_0}^{\R}(t_k)$ for different values of $h_0$. The next table shows the values of $h_0$ giving good results together with the value  $\rho_{h_0}^{\R}(t_{k'_2})$ at which the sequence became stable and the new maximum $\gamma'_2$.

\begin{center}
 Table 8:
\end{center}
\begin{center}
\begin{tabular}{|c|c|c|c|}
  \hline
$i$ & $h_0$ &  $\rho_{h_0}^{\R}(t_{k'_2})$ & $\gamma'_2$ \\
\hline
5   & $1$  & 0.05275311956 & 0.05275311956\\
11  & $1$    & 0.05120027918 & 0.05120027918\\
14  & $1$  & 0.06065146651 & 0.06065146651\\
15  & $1$  & 0.09855121223 & 0.09855121223\\
18  & $ 1$ & 0.08937284288 & 0.08937284288 \\
19  & $4/5$ & 0.1099404739  & 0.1099404739 \\
20  & $ 1$ & 0.07070538339 & 0.07070538339 \\
22  & $ 1$ & 0.06827542251 & 0.06827542251 \\
26  & $ 1$ & 0.07787356026 & 0.07787356026 \\
 \hline
\end{tabular}
\end{center}

\vspace{0.5cm}{\bf \noindent Acknowledgements.}
The authors deeply thank Sonia P\'erez-D\'{\i}az and J. Rafael Sendra for many useful discussions on
the topics treated in this paper.
\vspace{0.5cm}

{

\end{document}